\documentclass[a4paper,12pt]{amsart}
\usepackage{amssymb}
\usepackage{amsmath}
\usepackage{stmaryrd}
\usepackage{amscd,amsthm,amssymb}
\usepackage{enumerate}
\usepackage{color}

\scrollmode
\usepackage{latexsym}

\def\.{\cdot}

\def\vs{\vskip .6cm}

\def\beq{\begin{equation}}
\def\eeq{\end{equation}}
\def\bea{\begin{eqnarray*}}
\def\eea{\end{eqnarray*}}
\def\beaa{\begin{eqnarray}}
\def\eeaa{\end{eqnarray}}
\def\ba{\begin{array}}
\def\ea{\end{array}}

\def\id{\mathrm{id}}
\def\be{\begin{equation}}
\def\ee{\end{equation}}
\def\tr{\mathrm{tr}}
\def\Aut{\mathrm{Aut }}

\def\so{\mathfrak{so}}
\def\su{\mathfrak{su}}

\def\SU{\mathrm{SU}}
\def\U{\mathrm{U}}

\def\G{\mathrm{G}}

\def\End{\mathrm{End}}

\def\Spin{\mathrm{Spin}}

\def\scal{\mathrm{scal}}

\def\Id{\mathrm{id}}

\def\S{\mathrm{S}}

\newtheorem{epr}{Proposition}[section]
\newtheorem{ath}[epr]{Theorem}

\newtheorem{ecor}[epr]{Corollary}

\theoremstyle{definition}
\newtheorem{ede}[epr]{Definition}

\title[Conformal Killing forms]{Conformal Killing forms on nearly K\"ahler manifolds}

\author{Antonio M. Naveira, Uwe Semmelmann}
\address{
Antonio M. Naveira\\
Departamento de Matem\'aticas, Facultad de Ciencias Matem\'aticas, Universidad de Valencia E.G.\\
Av. Vicente Andr\'es Estell\'es 1\\
46100 Burjassot, Valencia, Spain}
\email{naveira@uv.es}

\address{Uwe Semmelmann\\
Institut f\"ur Geometrie und Topologie \\
Fachbereich Mathematik\\
Universit{\"a}t Stuttgart\\
Pfaffenwaldring 57 \\
70569 Stuttgart, Germany
}
\email{uwe.semmelmann@mathematik.uni-stuttgart.de}

\date{\today}
\begin{document}

\begin{abstract}
We study conformal Killing forms on compact $6$-dimensional nearly K\"ahler manifolds. Our main result
concerns forms of degree $3$. Here we
give a classification showing that all conformal Killing $3$-forms are linear combinations of  $d \omega$ 
and its Hodge dual  $\ast d \omega$, where $\omega$ is the fundamental $2$-form of the nearly
K\"ahler structure. The proof is based on a fundamental integrability condition for conformal
Killing forms. We have partial results in the case of conformal Killing $2$-forms. In particular we show the non-existence of
$J$-anti-invariant Killing $2$-forms.
\vs
\end{abstract}
\maketitle

\section{Introduction}

Conformal Killing forms are a natural generalization of conformal vector fields on Riemannian manifolds.
They are defined as sections in the kernel of a certain conformally invariant first order differential operator.
Conformal Killing forms which are in addition co-closed are called Killing forms. Equivalently Killing forms are characterized
by a totally skew-symmetric covariant derivative, thus generalizing Killing vector fields in 
degree one. Killing forms were intensively studied in physics since they define first integrals of the equations of motion, i.e. functions which are constant along
geodesics. 

There are not too many known examples of compact Riemannian manifolds admitting
non-parallel conformal Killing forms and even less for Killing forms.  On the standard sphere the space of conformal Killing forms 
is of maximal dimension. It coincides with the sum of two eigenspaces of the Laplace operator on forms corresponding to minimal 
eigenvalues.
Otherwise examples of Killing forms are usually related
to special geometric situations, e.g. Killing forms exist on  Sasakian, nearly K\"ahler or nearly parallel 
${\mathrm G}_2$ manifolds. 

It is known that there are no non-parallel Killing forms on compact manifolds
of special holonomy, e.g. for  symmetric spaces, K\"ahler manifolds,  or manifolds of holonomy ${\mathrm G}_2$ or $\Spin_7$
(cf. \cite{andrei-uwe5},  \cite{uwe2}, \cite{yama3}).
Much less is known about conformal Killing forms on such manifolds or  on manifolds with a weakened holonomy condition, i.e. with some special structure group reduction. There is a classification of conformal Killing forms on compact K\"ahler manifolds (cf. \cite{andrei-uwe4}).
In particular, in degree $2$ they turn out to be $J$-invariant and related to so-called Hamiltonian 
$2$-forms, e.g. studied in \cite{acg}. In the more general case there are only partial and
unpublished results, e.g. in \cite{diss} it could be shown that any conformal Killing form on a compact
Sasaki-Einstein manifold has to be one of the standard examples.

In the present article we study conformal Killing forms on compact $6$-dimensional nearly K\"ahler manifolds. These are by definition almost
Hermitian manifolds $(M,g,J)$ where the fundamental $2$-form $\omega$ is a  Killing $2$-form, i.e. has
totally skew-symmetric covariant derivative. In  dimension $6$ nearly K\"ahler manifolds are particularly interesting,
e.g. they are Einstein and the metric cone has holonomy ${\rm G}_2$. Moreover, with $\omega$ also $d\omega$ is
a conformal Killing form and, since the Hodge $*$-operator  preserves the space of conformal Killing forms, also $\ast \omega$ and $\ast d \omega$ are conformal Killing forms. Our main result
states that there are no other examples in degree $3$, i.e. we have

\begin{ath}\label{main}
Let $(M^6, g, J)$ be a compact $6$-dimensional strict nearly K\"ahler manifold. Then any
conformal Killing $3$-form on $M$ is a linear combination of $d\omega$ and $\ast d \omega$.
\end{ath}

It is interesting to note that the proof of this theorem is based on a fundamental integrability condition for
conformal Killing forms. This condition, which is related to the fact that conformal Killing forms are 
components of a parallel section in the so-called prolongation bundle, is known for a long time but so far
stayed with almost no application.

In degree $2$ the situation seems to be more complicated. Here it is not possible to apply the integrability condition
and also other methods useful in the case of degree $3$ can not be applied. However,  it is still possible to show the non-existence of Killing $2$-forms of
special type, e.g. we show in Propositions \ref{prop2} and  \ref{prop3}  that on a compact $6$-dimensional strict nearly
K\"ahler manifold there are no $J$-anti-invariant Killing  or $*$-Killing forms of degree $2$.

In a recent work  I. Dotti and C. Herrera study  Killing $2$-forms on homogeneous spaces $G/K$ invariant under the group $G$.
In particular they can show that there are no other Killing $2$-forms on the nearly K\"ahler  flag manifold $\SU(3)/T^2$ (cf.  \cite{dh}).
Given this result it seems to be a natural to  conjecture that compact strict nearly K\"ahler manifold in dimension $6$ do not admit
Killing $2$-forms non-proportional to the fundamental $2$-form $\omega$.

{\sc Acknowledgment.} We would like to thank Andrei Moroianu for his interest in our work and for helpful comments.
The first author has been partially supported by MINECO-FEDER grant MTM2016-77093-P and Generalitat Valenciana Project PROMETEO0II/2014/064.


\section{Nearly K\"ahler manifolds}

Let $(M,g,J)$ be an almost Hermitian manifold, i.e. the Riemannian manifold $(M,g)$ admits an  almost
complex structure $J$ compatible with $g$.  In this situation the fundamental $2$-form is defined as  $\omega(X,Y) = g(JX, Y) $
for tangent vectors $X, Y$.  
Using a local orthonormal basis $\{e_i\}$ it can be written as
$\,
\omega = \frac12 \sum e_i \wedge Je_i
$.

\begin{ede}[cf. \cite{gray}] 
A nearly K\"ahler manifold is an  almost Hermitian manifold $(M,g,J)$ such that the fundamental $2$-form
$\omega$ is a Killing $2$-form, i.e. satisfies $X \lrcorner \nabla_X \omega = 0$ for every tangent vector $X$. A nearly K\"ahler 
manifold is called strict if $\nabla_XJ \neq 0$ for every $X \neq 0$.
\end{ede}

In this article we are mainly interested in strict $6$-dimensional nearly K\"ahler manifolds. These manifolds are automatically Einstein with positive scalar curvature. As usual we will  normalize the scalar curvature to
$\scal = 30$. The $3$-form $\Psi^+ := d\omega$ is of type $(3,0) + (0,3)$ and the real part of a complex
volume form $\Psi := \Psi^+ + i \Psi^-$, where $\Psi^-$ is the Hodge dual of $\Psi^+$. The differential forms
$\omega$ and $ \Psi^+$ define a structure group reduction from ${\mathrm O}(6)$ to $\SU(3)$. Conversely, $6$-dimensional 
nearly K\"ahler manifolds can be described as Riemannian manifolds $(M^6,g)$ with a $\SU(3)$-structure
$(\omega, \Psi^{\pm})$ satisfying the equations
$
d \omega = 3 \Psi^+
$
and
$
d \Psi^-  = - 2 \, \omega \wedge \omega
$
.
More details on SU(3)-structures and 6-dimensional nearly K\"ahler manifolds can be found in \cite{andrei-uwe3}.

In dimension $6$ there are only a few examples known. 
The homogeneous ones, classified in \cite{butruille},
are  the standard sphere $\S^6$ and $\S^3 \times \S^3, \SU(3)/T^2$ and ${\mathbb C}P^3$ with their nearly K\"ahler metrics. Only quite recently new examples were found. These are metrics of cohomogeneity one on $\S^6$ and on $\S^3 \times \S^3$ (cf. \cite{fh18}).

In the rest of this section we collect some results on the curvature of nearly K\"ahler manifolds, which are needed later.
Let $(M^n, g, J)$ be a nearly K\"ahler manifold. The Riemannian curvature $R$ can  be
considered as a map $R :\Lambda^2TM \rightarrow \Lambda^2TM$, the so-called curvature operator. It is defined by
the equation $ g(R (X \wedge Y), Z\wedge V ) = R(X,Y,Z, V)$ for any tangent vectors
$X, Y, Z, V$. Note that with this convention $R$ acts as $-\Id_{\Lambda^2T}$ on the standard sphere $S^n$.

For further use  we  have to recall the definition  of the curvature endomorphism $q(R)$.
Let 
$P=P_{SO(n)}$ be the frame bundle and $EM$ any vector bundle associated to $P$ via a
$SO(n)$-representation $\rho: SO(n) \rightarrow \Aut (E)$, where $E$ is some real or complex vector space.
Then $q(R) \in \End (EM)$ is defined as
\be\label{qR}
q(R) \;:=\; \tfrac12 \,  (e_i \wedge e_j)_\ast \circ  R(e_i \wedge e_j)_\ast \ ,
\ee
where $\{e_i\}, i = 1, \ldots , n$, is a local orthonormal frame. 
Here and henceforth, we use Einstein’s summation convention on repeated subscripts. 
We also identify $TM$ with $T^*M$ using the metric.
The $2$-form $X\wedge Y \in \Lambda^2 T \cong \so(n)$
acts via the differential $\rho_\ast$ of the representation $\rho$,  we write  $(X \wedge Y)_\ast = \rho_\ast (X \wedge Y)$. In particular we have
for any tangent vectors $X,Y, Z$ the standard action of $\so(n)$ written as
$
(X\wedge Y)_\ast\, Z = g(X,\,Z) \,Y - g ( Y,\, Z) \, X
$.
Moreover for any section $\varphi \in \Gamma(EM)$ we have
$
R (X\wedge Y)_\ast \, \varphi \;=\; R_{X,Y}\,\varphi
$.
It is easy to check that $q(R)$ acts as  the Ricci endomorphism on tangent vectors. We remark that
$q(R)$ may be defined in this way for any curvature tensor $R$, e.g. for the curvature tensor
$\bar R$ of the canonical hermitian connection $\bar \nabla$.
Recall that the canonical connection $\bar \nabla$ is defined by
$
\bar \nabla_X Y := \nabla_XY - \tfrac12 J(\nabla_XJ)Y
$
for any vector fields $X, Y$. This is a $\U_m$-connection, i.e. $\bar \nabla g = 0$ and $\bar \nabla J = 0$.

It is well-known that the Riemannian curvature tensor $R$ of a nearly K\"ahler manifold can be
written as $R_{X\,Y} = - (X \wedge Y)_* \,+\, R^{CY}_{X\,Y} $,
where $R^{CY}_{X\,Y} $ is a 
curvature tensor of Calabi-Yau type (cf. \cite{andrei-uwe1}, p. 253). In other words, the
curvature operator $R^{CY}= R + \id : \Lambda^2T \rightarrow \Lambda^2T$ takes
values only in $\Lambda^{1,1}_0T \cong \su(n)$. Equivalently we have for  every tangent 
vectors $X, Y$ the equations
\beq\label{curvature}
R^{CY}_{X\,Y} \circ J \;=\; J \circ R^{CY}_{X\,Y}
\qquad\mbox{and}\qquad
\tr \, ( R^{CY}_{X\,Y} \circ J) \;=\; 0 \ .
\eeq
Note, that the second equation can be written as $R^{CY}(\omega) =0$ and  also as $R(\omega) = -\omega$.
The first equation in \eqref{curvature} can be rewritten in several ways, giving compatibility equations for the
Riemannian curvature $R$ and the almost complex structure $J$. In particular we have

\begin{epr}\label{comp}
Let $(M^n, g, J)$ be a nearly K\"ahler manifold. 
Then for any tangent vectors $X, Y, Z$ the following equations hold
\begin{enumerate}
\item[(i)]
\quad
$ R_{X\,Y} JZ \;=\; J \, R_{X\,Y} Z \;+\; J \, (X \wedge Y)_\ast \,Z \; - \; (X\wedge Y)_\ast \,JZ $
\medskip
\item[(ii)]
\quad
$
R_{JX\, JY} Z \;=\; R_{X\,Y}Z \;+\; (X\wedge Y)_\ast \, Z \;-\; (JX\wedge JY)_\ast \,Z
$
\medskip
\item[(iii)]
\label{comp1}
\quad
$
R_{X\, JY} Z \;=\; -\,R_{JX\,Y}Z \;-\; (JX\wedge Y)_\ast \,Z \;-\; (X \wedge JY)_\ast \, Z 
$
\end{enumerate}
\end{epr}
\proof
The first equation is exactly \eqref{curvature} written in terms of $R$. The second equation
follows from the first by the following calculation
\bea
g(R_{JX\,JY} Z, V) &=& g(R_{Z \, V} JX, JY)
\;=\;
g(R_{Z\, V}X, Y) \;+\; g( J (Z \wedge V)_\ast \, X \,-\, (Z \wedge V)_\ast \, JX, JY) \\[1ex]
&=&
g(R_{X\,Y}Z, V ) \;+\; g((Z\wedge V)_\ast \, X, Y) \;-\; g((Z\wedge V)_\ast\,JX,  JY)\\[1ex]
&=&
g(R_{X\,Y}Z, V ) \;+\; g(Z\wedge V, X\wedge Y) \;-\; g(Z\wedge V, JX\wedge JY )\\[1ex]
&=&
g(R_{X\,Y}Z, V ) \;+\; g((X\wedge Y)_\ast Z, V) \;-\; g((JX \wedge JY)_\ast Z, V) \ .
\eea
The third equation follows from the second by replacing $X$ with $JX$ \qed

\medskip

As a consequence of these curvature equations we have several equations on the space of $p$-forms, which
we will later apply to conformal Killing forms on nearly K\"ahler manifolds. For any $p$-form $\sigma$ 
and any tangent vector $X$ we define the curvature expressions
\be\label{curvex}
R^+(X) \,\sigma\;:=\;  e_i \wedge R_{X\, e_i}\,\sigma
\qquad \mbox{and}\qquad
R^-(X) \,\sigma\;:=\;   e_i \, \lrcorner \, R_{X\,e_i}\,\sigma \ ,
\ee
where $\{ e_i\}$ is a local orthonormal frame. 
Recall that the almost complex structure $J$ extends to a map on $p$-forms. It is defined as
$
J_*\,\sigma =  Je_i \wedge e_i \, \lrcorner \,\sigma = \omega_\ast \, \sigma
$.
In particular, we note that $R(\omega) = - \omega$ can be written as $R_{e_i \, Je_i}\sigma  \;=\; -\,2\, J_*\,\sigma$.
Moreover, we recall that $(J_*)^2$ acts as $-(p-q)^2 \,\id$ on forms of type $(p,q) + (q,p)$ (cf. \cite{andrei-uwe2}, p. 60).
For later use we note that $J_*$ is injective on the space of $3$-forms. Indeed, 
 the space of $3$-forms on an almost Hermitian manifold decomposes into the sum of spaces of forms of types $(3,0)+(0,3)$ and  $(1,2) +(2,1)$.

\begin{ecor}\label{Folgerung}
Let $(M^6, g, J)$ be a $6$-dimensional nearly K\"ahler manifold. 
Then for any tangent vector $X$ and any $p$-form $\sigma$ the following equations hold 
\begin{enumerate}
\item[(i)]
$
\quad
 Je_i \,\lrcorner\, R_{X\, e_i} \,\sigma \;=\; R^-(JX)\,\sigma \;-\; J_*\, (X \,\lrcorner  \,\sigma ) \;+\; 2\, X \wedge \omega\,\lrcorner \,\sigma\;+\; (5-p)\,JX\,\lrcorner\,\sigma
$
\medskip
\item[(ii)]
$
\quad
Je_i \lrcorner \, R^+(e_i)\,\sigma \;+\;  Je_i \wedge R^-(e_i) \,\sigma  \;=\;-\, 6\, J_*\,\sigma
$\label{Kern}
\end{enumerate}
\end{ecor}
\proof
We use the third equation  of Proposition \ref{comp} and the definition of $J$  to obtain
\bea
Je_i \,\lrcorner\, R_{X, \,e_i}\,\sigma &=& - e_i \,\lrcorner\,R_{X, \,Je_i} \,\sigma\\[1ex]
&=&
\quad R^-(JX) \,\sigma \;+\;  e_i \,\lrcorner\, ( e_i \wedge JX\,\lrcorner \,  \;-\; JX \wedge e_i \lrcorner \,)\,\sigma\\[0ex]
&&  \phantom{xxxxxxxxxxxxxxx}
+\;  e_i \,\lrcorner\, (Je_i \wedge X \lrcorner \;-\; X \wedge Je_i \lrcorner \, )\,\sigma  \\[1ex]
&=&
R^-(JX)\,\sigma \;+\;    6 \, JX \,\lrcorner \,\sigma \,-\, (p-1)\,JX\,\lrcorner \,\sigma \;-\; JX \,\lrcorner \,\sigma   \\
&& \phantom{xxxxxxxxxxxxx}
-\; J(X\,\lrcorner \,\sigma)   \;-\; JX \,\lrcorner\,\sigma \;+\; X \wedge 2\,\omega \,\lrcorner\,\sigma
\\[1ex]
&=&
R^-(JX)\,\sigma \;-\; J(X\,\lrcorner \,\sigma) \;+\; 2\,X \wedge \omega \,\lrcorner\,\sigma \;+\; (5-p)\,JX\,\lrcorner\,\sigma
\eea
This proves the first equation of Corollary \ref{Folgerung}. Rewriting the first summand of the second equation we find
 $$
 Je_i \, \lrcorner \, R^+(e_i)\,\sigma 
 \;=\;
 Je_i \,\lrcorner \, ( e_j \wedge R_{e_i \, e_j})\,  \sigma   \\[-.8ex]
 \;=\;
g(Je_i, \, e_j) \, R_{e_i \, e_j}\, \sigma \,+\,  e_j \wedge Je_i \, \lrcorner \, R_{e_j \, e_i} \sigma   \ .\\[1ex]
$$
Using the first equation of Corollary \ref{Folgerung} with $X=e_j$ and the equation  $ R_{e_i Je_i} = -2J $  we obtain
 \bea
 Je_i \, \lrcorner \, R^+(e_i)\,\sigma 
 &=&
  R_{e_i \, Je_i} \, \sigma 
 \;+\;
  e_j \wedge \left(R^-(Je_j)\sigma  \,-\, J(e_j \,\lrcorner \, \sigma) \,+\, 2 \, e_j \wedge \omega \,\lrcorner \, \sigma
      \,+\, (5-p)Je_j \, \lrcorner \sigma   \right)\\[1ex]
&=&
-\,2\,J \,\sigma \;-\; Je_j \wedge R^-(e_j) \,\sigma 
\;-\; e_j \wedge J(e_j \, \lrcorner \,\sigma) \;-\; (5-p)\,Je_j \wedge e_j \,\lrcorner \, \sigma\\[1ex]
&=&
-\, Je_j \wedge R^-(e_j) \,\sigma \;-\; (7-p)\, J\, \sigma \;-\; e_j \wedge J(e_j \, \lrcorner \,\sigma) 
 \eea 
It remains to determine the sum 
$
\; e_j \wedge J(e_j \, \lrcorner \,\sigma) =  e_j \wedge Je_i \wedge e_i \, \lrcorner \,e_j \,\lrcorner \, \sigma
$.
Here we compute
\bea
e_j \wedge J(e_j \, \lrcorner \,\sigma) 
&=&
-\, Je_i \wedge e_j \wedge e_i \,\lrcorner \, e_j \,\lrcorner \, \sigma
\;=\;
 Je_i \wedge e_i \,\lrcorner \, (e_j \wedge e_j \,\lrcorner \,\sigma) \;-\; Je_i\wedge e_i \, \lrcorner \, \sigma\\[1ex]
&=&
(p-1)\, J\,\sigma 
\eea
Substituting this result into the last expression for $Je_i \, \lrcorner \, R^+(e_i)\,\sigma $ we obtain the second
equation of Corollary \ref{Folgerung}.
\qed


\section{Conformal Killing forms}

In this section we will recall the definition, examples and important properties of conformal Killing forms. More details can
be found in \cite{uwe}.

Conformal Killing $p$-forms on a Riemannian manifold $(M^n,g)$ are defined as differential forms $\sigma \in \Gamma(\Lambda^p TM)$ satisfying
for any tangent vector $X$  the equation
$$
\nabla_X \,\sigma \;=\; \tfrac{1}{p+1} \, X \lrcorner \, d \,\sigma \;-\; \tfrac{1}{n-p+1}\, X \,\wedge \, d^* \,\sigma \ .
$$

Conformal Killing forms which in addition are co-closed are called Killing forms. Closed conformal
Killing forms are also called $\ast$-Killing forms. Indeed the Hodge $\ast$-operator preserves the space
of conformal Killing forms and maps Killing forms to $\ast$-Killing forms and vice versa. 

Every parallel form is trivially a Killing form. For $p=1$ Killing forms are dual to Killing vector fields. 
The standard sphere $(\S^n, g_0) $ is the compact manifold with the maximal number of conformal Killing forms. Here every
conformal Killing form is a linear combination of a Killing and a $\ast$-Killing form. The space of Killing
forms on $\S^n$ coincides with the eigenspace of the Laplace operator on co-closed forms for the minimal eigenvalue.
Other interesting examples are related to special geometric structures. For nearly K\"ahler manifolds the fundamental
$2$-form is by definition a Killing form. Similarly, the defining $3$-form of a nearly parallel $\G_2$-structure is by
definition a Killing form. If $\eta$ is the $1$-form dual to the Reeb vector field $\xi$ defining a Sasakian structure then
all the forms $\eta \wedge d\eta^k$ for $k=0, \ldots , n$ are Killing forms.

These examples on the standard sphere,  $6$-dimensional nearly K\"ahler,  nearly parallel $\G_2$- and on 
Sasakian manifolds are so-called special Killing forms. They are Killing forms $\sigma$ satisfying  the additional equation
$
\nabla_X d \, \sigma = c \,X \wedge \sigma
$
for some real constant $c$ and every tangent vector $X$. 
Special Killing forms on compact Riemannian manifolds $(M^n, g)$ 
were classified in \cite{uwe}. They turn out to be in  bijective correspondence to parallel forms on the metric cone over $M$.
In particular it follows that on $6$-dimensional nearly K\"ahler manifolds
$d\omega$ is a closed conformal Killing form and that $\Delta \omega = 12 \omega$ 
(cf. \cite{uwe}, Prop. 4.2).
As far as we know, the only examples of non-parallel, non-special Killing forms in degree larger than one are the 
fundamental $2$-forms of  strict nearly K\"ahler manifolds  in dimension larger than $6$ and the torsion forms of 
metric connections with skew-symmetric and parallel torsion, e.g. the torsion form on naturally reductive spaces.

In this article we will consider conformal Killing forms on compact manifolds. Here one has a rather useful additional
characterization.  First, we recall from \cite{uwe} that  conformal Killing $p$-forms $\sigma$  on an arbitrary Riemannian 
manifold $(M^n, g)$ satisfy a second  order Weitzenb\"ock equation:
\be\label{qr}
q(R) \, \sigma
\;=\;
\tfrac{p}{p+1} \, d^* d \, \sigma \;+\; \tfrac{n-p}{n-p+1}\, d d^\ast \, \sigma \ ,
\ee
where $q(R) \in \End (\Lambda^pTM)$ is the  curvature term defined in \eqref{qR}.

Let $\Delta = d^* d + d d^\ast$ denote the Hodge-Laplace operator on
forms. Then  it follows from \eqref{qr} that $\Delta \,\sigma = \frac{p+1}{p} q(R) \,\sigma$  holds for Killing $p$-forms $\sigma$
and $\Delta\, \sigma = \frac{n-p+1}{n-p} q(R) \,\sigma$ for  $\ast$-Killing $p$-forms. A simple
integration argument shows that these equations characterize Killing resp. $\ast$-Killing forms
on compact manifolds, i.e.  a $p$-form $\sigma$ is a Killing form if and only if $d^*\sigma = 0$ and
 $\Delta \,\sigma = \frac{p+1}{p} q(R) \sigma$. Similarly $\sigma$ is a $\ast$-Killing form
 if and only if $d\sigma =0$ and $\Delta \sigma = \frac{n-p+1}{n-p} q(R) \sigma$.
Moreover, we have a similar characterization for conformal Killing forms in middle dimension, i.e.
 a $m$-form $\sigma$ on a $2m$-dimensional compact manifold $M$ is
conformal Killing if and only if the equation $\Delta \sigma = \frac{m+1}{m}q(R)\sigma$ is satisfied  (cf.  \cite{uwe}, Cor. 2.5).

Later we will study conformal Killing forms on $6$-dimensional nearly K\"ahler manifolds, which are automatically Einstein. 
On Einstein manifolds more can be said about conformal Killing forms, e.g. it is known that in this situation for any
conformal Killing $2$-form $\sigma$ the codifferential $d^* \sigma$ is dual to a Killing vector field.
However, here we are more interested in a commutator rule for the differential and codifferential with the
curvature endomorphism $q(R)$. We have the following result.
Let $(M^n, g)$ be an Einstein manifold. Then every conformal Killing $p$-form $\sigma$ satisfies  the
two equations 
\beq\label{Einstein}
d(q(R)\,\sigma) \;=\; \tfrac{scal}{n} \, d \, \sigma \;+\; \tfrac{p-1}{p+1}\, q(R)\, d\,\sigma  
\eeq
and
\beq\label{Einstein2}
\quad d^*(q(R)\,\sigma) \;=\;  \tfrac{scal}{n} \, d^* \, \sigma \;+\; \tfrac{n-p-1}{n-p+1}\, q(R)\, d^*\,\sigma \ .
\eeq
The statements follow from easy local calculations (e.g. cf. \cite{uwe1}, Prop. 4.4.12 and Cor. 7.1.2).

Finally we want to mention an important integrability condition for conformal Killing forms. Indeed, the conformal Killing
equation is of finite type, i.e. there is a finite prolongation and conformal Killing forms are components of a parallel 
section in a larger bundle, with respect to a suitable connection. In particular, the curvature of this connection
vanishes on conformal Killing forms and one component of the corresponding equation is the following 
integrability condition (cf.  \cite{uwe}, Prop. 6.4). Let $\sigma$ be a  conformal Killing $p$-form then
\beaa \label{int}
R_{X\,Y} \sigma &=&
\tfrac{1}{p(n-p)} \, (X\wedge Y)_\ast \,q(R) \,\sigma
\;+\;
\tfrac{1}{p}\, (\, Y\, \lrcorner\, R^+(X ) \,\sigma \;-\; X \,\lrcorner \, R^+(Y)\,\sigma \,) \\[1ex]
&&
\phantom{xxxxxxxxxxxxxxxxxxxxxxx}
\;+\;
\tfrac{1}{n-p}\,(\, Y \wedge R^-(X)\,\sigma \;-\; X \wedge R^-(Y)\,\sigma  \, ) \nonumber
\eeaa
is satisfied for all  tangent vectors $X, Y$. The curvature expressions $R^+(X)$ and $R^-(X)$ were defined in
\eqref{curvex}. We will use this equation for conformal Killing $3$-forms on 
$6$-dimensional nearly K\"ahler manifolds.

\bigskip


\section{Conformal Killing $3$-forms
on nearly K\"ahler manifolds}

Let $(M^6, g, J)$ be a  compact $6$-dimensional strict nearly K\"ahler manifold with scalar curvature 
normalized to $\scal = 30$ and let $\sigma$ be a conformal
Killing $3$-form on $M$. In this section we will give the proof of our main theorem, i.e. Theorem \ref{main}, and show that
that the $3$-form $\sigma $ has to be a linear combination of  $d \, \omega$ and $\ast d \,\omega$.

In the first step we use the integrability condition   \eqref{int} with $X=e_i, Y= J e_i$, for an orthonormal
frame $\{e_i\}$. After summation we obtain
\beq \label{3curvature}
 R_{e_i \, Je_i} \sigma \;=\; \tfrac19 \, (e_i \wedge Je_i)_\ast \, q(R) \, \sigma
\;+\; \tfrac23\, Je_i \,\lrcorner \, R^+(e_i)\,\sigma \;+\; \tfrac23 \, Je_i \wedge R^-(e_i)\, \sigma \ .
\eeq
Since  $R(\omega) = - \omega$ for the fundamental $2$-form $\omega$, the left-hand side of this equation can be rewritten as 
$
R_{e_i \, Je_i} \sigma = 2 \, R(\omega) \sigma = -2 \, \omega_\ast \sigma =
- 2 \, J \sigma $.
The first summand on the right-hand side
is equal to $\frac29 J q(R) \sigma$. Moreover, substituting the second equation of  Corollary \ref{Folgerung} into  \eqref{3curvature}, 
we see that \eqref{3curvature} is equivalent to the equation
$
J q(R) \sigma = 9 J \sigma
$.
But $J$ is injective on $3$-forms, thus we find $q(R) \sigma = 9 \sigma$ for every conformal Killing $3$-form $\sigma$.

\medskip

Let $\sigma$ be a conformal Killing $3$-form on a $6$-dimensional nearly K\"ahler manifold. Then $\sigma$ is a conformal Killing form 
in middle dimension and  we have
$
\Delta \sigma = \frac43 \, q(R)\, \sigma = 12 \, \sigma
$.
Since the Laplace operator $\Delta = d \, d^* + d^* d$ commutes with $d$ and $d^*$,  also
$\Delta d\sigma = 12 d \sigma$, \;$\Delta d^*\sigma = 12 d^*\sigma$ and $\Delta d \, d^*\sigma = 12 d\, d^*\sigma$ hold.
Using  \eqref{Einstein} and \eqref{Einstein2} we obtain 
\beq
d (q(R) \sigma) \;= \;
5 \, d\, \sigma \;+\; 
\tfrac12\, q(R)\, d\,\sigma 
\qquad \mbox{and} \qquad
d^* (q(R) \sigma) \;= \;
5 \, d^*\, \sigma \;+\; 
\tfrac12\, q(R)\, d^*\,\sigma  \ .
\eeq
Since   $q(R) \sigma = 9\, \sigma$  these equations imply
$$
q(R)\, d\sigma \;=\; 8\, d\sigma
\qquad \mbox{and} \qquad
q(R) \, d^\ast \sigma \;=\; 8\, d^\ast  \sigma \ .
$$
and it follows that
$
\Delta\,  d^*\sigma = 12 \, d^*\sigma = \frac32 \, q(R) \, d^*\sigma
$
and similarly that
$
\Delta \,d \sigma =  \frac32 \, q(R) \, d\sigma
$.
Hence, by the characterization of Killing and $\ast$-Killing forms on compact manifolds given above,
we conclude that for every conformal Killing $3$-form $\sigma$ the form 
$d^*\sigma$ is a Killing $2$-form and $d\sigma$ is a $\ast$-Killing $4$-form.

Now, let  $\sigma$ be a Killing form, i.e. a co-closed conformal Killing  $3$-form, then $d\sigma$ is a $\ast$-Killing $4$-form 
and the $\ast$-Killing equation for $d \sigma$ reads for any tangent vector $X$ as
$$
\nabla_X d \sigma \;=\; - \, \tfrac13\, X\,\wedge \, d^* \, d \, \sigma \;=\; -\, \tfrac13\, X\wedge\,  \Delta\,  \sigma
\;=\; - \, 4\, X\, \wedge \, \sigma
$$
and similarly 
$
\nabla_Xd^*\sigma = 4\, X\, \lrcorner \, \sigma
$.
Thus we see that every Killing $3$-form $\sigma$  already is a
special Killing $3$-form and similarly every  $\ast$-Killing $3$-form   
is the Hodge dual of a special Killing $3$-form.

%

In general a conformal Killing $3$-form $\sigma$ needs not to be closed or co-closed. However,
it follows, that the $3$-form $\sigma - \frac{1}{12}dd^*\sigma$ is co-closed and 
$\sigma - \frac{1}{12} d^* d \, \sigma$ is closed. Indeed
$$
d^*(\sigma \,-\, \tfrac{1}{12} \, d \,d^* \sigma) \;=\; d^* \sigma - \tfrac{1}{12} \,d^* d \,d^* \sigma \;=\; d^*\sigma
\;-\; \tfrac{1}{12} \,\Delta \, d^* \sigma \;=\; 0 \ .
$$
A similar calculation shows that $\sigma - \frac{1}{12} d^* d \, \sigma$ is closed. Using \eqref{Einstein} once again, this time for the Killing $2$-form $d^\ast \sigma$, we obtain
$
q(R) d d^\ast \sigma \;=\; 9 \, d d^\ast \sigma
$.
Hence, $ \Delta \, d d^\ast \sigma  = 12\,  \sigma =  \frac43\, q(R) d d^\ast \sigma$ and we see 
as above that $d d^\ast \sigma$ has to be a closed conformal Killing $3$-form. But then 
$\sigma - \frac{1}{12}dd^*\sigma$ is a conformal Killing form too and thus, since co-closed, 
it is in fact a Killing $3$-form. 
%
%
We showed already that Killing $3$-forms are automatically special. Thus we conclude that 
$\sigma - \frac{1}{12}dd^*\sigma$  is a special Killing $3$-form. In the same way we can show
that  $\sigma - \frac{1}{12} d^* d \, \sigma$ is  the Hodge dual of a special Killing $3$-form.

Special Killing forms were classified in \cite{uwe}, Thm. 4.8. In particular,   the only special Killing $3$-forms on
$6$-dimensional nearly K\"ahler manifolds are constant multiples of $\ast d \omega$. It follows that
$
 \sigma \,-\, \tfrac{1}{12} \, d d^* \sigma = \lambda \, \ast d \omega
$
for some real constant $\lambda$. Similarly we have that  $  \sigma \,-\, \tfrac{1}{12} \, d^\ast d \sigma$
is the Hodge dual of a special Killing $3$-form and thus
$
 \sigma \,-\, \tfrac{1}{12} \, d^* d \sigma = \mu  \, d \omega
$
for some real constant $\mu$. Since $\Delta \sigma = 12 \sigma$  we can write $\sigma$ as
$$
 \sigma
 \;=\;
  2\,  \sigma \;-\; \tfrac{1}{12} \, \Delta \sigma
 \;=\;
 (\sigma \,-\, \tfrac{1}{12} \, d^* d \sigma)
 \;+\;
  (\sigma \,-\, \tfrac{1}{12} \, d d^* \sigma)
  \;=\;
 \lambda  \ast d \omega  \,+\, \mu  \, d \omega
$$ 
This finishes the proof of Theorem \ref{main}.

%


\bigskip

\section{Killing 2-forms on nearly K\"ahler manifolds}

In this last section we want to make a few remarks concerning conformal Killing $2$-forms on nearly K\"ahler
manifolds. Here we can not use an argument similar to the one for $3$-forms. The coefficients of the last two
summands in our integrability condition \eqref{3curvature} are now different, so it is not possible to substitute 
the second equation of Corollary \ref{Folgerung}. Also, the extension of the almost complex structure $J$ 
is not injective on the space of all $2$-forms, it vanishes on $2$-forms of type $(1,1)$, i.e. $J$-invariant
$2$-forms.

However, it is still possible to exclude the existence of Killing and $*$-Killing $2$-forms of a special type, e.g. of type $(2,0)+(0,2)$,
i.e. anti-invariant $2$-forms. Let $\sigma$ be an anti-invariant $2$-form on a compact nearly K\"ahler manifold
$(M^6,g,J)$. Then there exists a vector field $\xi$ with $\sigma = \xi \,\lrcorner \, \Psi^+ = \Psi^+_\xi$. First we
compute the action of $q(R)$ on anti-invariant $2$-forms. Here we have 
$
q(\bar R)\, \Psi^+_\xi = \Psi^+_{q(\bar R)  \xi} =  4 \, \Psi^+_\xi 
$
(cf. \cite{andrei-uwe1}, p. 254 and Lemma 4.8). From Proposition 3.4 in \cite{andrei-uwe3} we then obtain
$
q(R) \,  \sigma = q(\bar R)\, \sigma + 4 \, \sigma = 8 \, \sigma 
$.

Assume $\sigma$ to be a Killing $2$-form then  we have
$
\Delta\,  \sigma = \tfrac32 \,q(R)\, \sigma = 12 \, \sigma
$.
But now we can use \eqref{Einstein} for the $3$-form $d\sigma$ to get
$
q(R) d\sigma = 9 d\sigma 
$
and thus
$
\Delta d \sigma = 12 \, \sigma = \tfrac43 q(R) d\sigma 
$.
Hence, $d\sigma$ is a closed conformal Killing $3$-form and by Theorem \ref{main} it has to be a 
constant multiple of $d\omega$, i.e. $d \sigma = \lambda d\omega$ for some real number $\lambda$.
Finally we note that $\Delta \sigma = d^* d \sigma = \lambda d^* d \omega = \lambda \Delta \omega = 12 \, \omega$.
This would imply $\omega = \sigma $, which is of course a contradiction. Summarizing we proved the
following

\begin{epr}\label{prop2}
Let $(M^6, g, J)$ be a compact nearly K\"ahler manifold. Then $M$ admits
no non-trivial $J$-anti-invariant Killing $2$-form.
\end{epr}

Finally we consider the case of closed $J$-anti-invariant  conformal Killing $2$-forms, i.e.
let  $\sigma := \Psi^+_\xi $ be a $\ast$-Killing $2$-form. Then by definition 
$
\nabla_X \sigma = - \tfrac13 \, X \wedge d^*\sigma
$
and in particular $\sigma$ is closed.
The characterization of  $\ast$-Killing $2$-forms implies
$
\Delta \sigma = \tfrac54 q(R) \sigma = 10 \sigma
$.
For the vector field  $V=d^*\sigma$ we have $d V =  d \, d^* \sigma = \Delta \sigma = 10 \sigma$. Hence,
$
\nabla_X dV = 10 \, \nabla_X \sigma = - \tfrac{10}{3} X \wedge V
$
and we obtain
$
\Delta V = d^* d V  = \tfrac{50}{3} V
$.
But with $\Delta \sigma = 10 \sigma$ we also have $\Delta V = 10 V$, which is  possible only with $V=0$.
Because $d V = 10 \sigma$ we conclude that also the conformal Killing $2$-form $\sigma$ has to vanish.
Summarizing we just proved

\begin{epr}\label{prop3}
Let $(M^6, g, J)$ be a compact nearly K\"ahler manifold. Then $M$ admits
no non-trivial $J$-anti-invariant $*$-Killing $2$-form.
\end{epr}

As already mentioned in the introduction, it follows from the work of  I. Dotti and C. Herrera in \cite{dh} that besides the fundamental
$2$-form there are no $\SU(3)$-invariant Killing $2$-forms on the nearly K\"ahler flag manifold
$\SU(3)/T^2$. It seems, to be a reasonable to conjecture that a Killing $2$-forms on a compact $6$-dimensional
strict nearly K\"ahler manifold has to be a skalar multiple of the fundamental $2$-form.

\end{document}